\documentclass[11pt]{article}

\usepackage{amsfonts,amsmath,amsthm,amssymb}
\usepackage{latexsym}

\newtheorem{theorem}{Theorem}[section]
\newtheorem{corollary}[theorem]{Corollary}
\newtheorem{lemma}[theorem]{Lemma}

\newtheorem{fact}[theorem]{Fact}

\theoremstyle{definition}

\newcommand{\beql}[1]{\begin{equation}\label{#1}}
\newcommand{\eeq}{\end{equation}}
\newcommand{\comment}[1]{}

\newcommand{\Abs}[1]{{\left|{#1}\right|}}

\newcommand{\Mean}{{\bf E}}

\newcommand{\Set}[1]{{\left\{{#1}\right\}}}

\newcommand{\RR}{{\mathbb R}}

\newcommand{\ZZ}{{\mathbb Z}}

\newcommand{\one}{{\bf 1}}
\newcommand{\zero}{{\bf 0}}

\newcommand{\inner}[2]{{\langle #1, #2 \rangle}}

\newcommand{\ft}[1]{\widehat{#1}}

\newcommand{\Li}{{\rm Li\,}}

\newcommand{\pty}{\mathbf {\oplus}}
\newcommand{\ptyc}{\bf \overline{\oplus}}
\newcommand{\eat}[1]{}

\newcounter{rem}
\newcounter{othm}
\def\theothm{\Alph{othm}} 
\newenvironment{othm}{
  \sf
  \vskip 0.10in
  \refstepcounter{othm}
  \noindent{\bf Theorem\ \theothm}
}{\vskip 0.10in}

\begin{document}

\title{Learning symmetric $k$-juntas in time $n^{o(k)}$}

\author{
{\sc Mihail N. Kolountzakis}\thanks{School of Mathematics, Georgia Institute of Technology, Atlanta GA 30332, USA, and
Department of Mathematics, Univ.\ of Crete, GR-71409 Iraklio, Greece. E-mail: {\tt kolount@gmail.com}.
Partially supported by European Commission IHP Network HARP
(Harmonic Analysis and Related Problems), Contract Number: HPRN-CT-2001-00273 - HARP
}
\and 
{\sc Evangelos Markakis}\thanks{Georgia Institute of Technology, Atlanta GA 30332, USA,
E-mail: {\tt \{vangelis, aranyak\}@cc.gatech.edu}} 
\and 
{\sc Aranyak Mehta}\footnotemark[2] 
 }

\date{April 2005}
\maketitle

\begin{abstract}
We give an algorithm for learning symmetric $k$-juntas
(boolean functions of $n$ boolean variables which depend only on an unknown set of $k$ of these variables) 
in the PAC model under the uniform distribution, which runs in time $n^{O(k/\log k)}$.
Our bound is obtained by proving the following result: Every symmetric boolean function on $k$ variables,
except for the parity and the constant functions, has a non-zero Fourier coefficient of order at least 1
and at most $O(k \bigl / \log k)$. This improves the previously best known bound of $3k \bigl / 31$~\cite{LMMV},
and provides the first $n^{o(k)}$ time algorithm for learning symmetric juntas.
\end{abstract}

\section{Introduction}
We consider a fundamental problem in computational learning theory:
learning in the presence of irrelevant information.  One formalization
of the problem is as follows: We want to learn an unknown boolean
function of $n$ variables, which depends only on $k \ll n$ variables
(typically $k$ is $O(\log{n})$). We call such a function a
$k$-junta. We are provided with a set of labelled examples
$\inner{x}{f(x)}$, where the $x$'s are picked uniformly and
independently at random from the domain $\{0,1\}^n$ (this is the PAC
model with uniform distribution). We wish to identify the $k$ relevant
variables and the truth table of the function.

The problem was first posed by Blum~\cite{Blum} and Blum and
Langley~\cite{BL}, and it is considered  \cite{blum-colt,MOS} to be one of the
most important open problems in the theory of uniform distribution
learning. It has connections with learning DNF formulas and
decision trees of super-constant size, see 
\cite{jun-bjt,jun-jackson,jun-mansour,jun-ver1,jun-ver2} for details. The general case is believed to be hard and has even been used to propose a cryptosystem \cite{bfkl}.
A trivial algorithm runs in time roughly $n^k$ by doing an exhaustive
search over all possible sets of relevant variables. Two important
classes of juntas are learnable in polynomial time: parity and
monotone functions. Learning parity functions can be reduced to
solving a system of linear equations over $\mathbb{F}_2$~\cite{HSW}.
Monotone functions have non-zero singleton Fourier coefficients (e.g.,
see~\cite{MOS}). For the general case, the first significant
breakthrough was given in~\cite{MOS} - learning with confidence
$1-\delta$ in time $n^{0.7k}{\rm poly}(2^k,n,\log{1/\delta})$.  Note
that we allow the running time to be polynomial in $2^k$, since this
is the size of the truth-table which is output.  In the typical
setting of $k=O(\log{n})$, this becomes polynomial in $n$.

In this paper we consider the class of {\em symmetric} $k$-juntas,
functions which are symmetric on their relevant variables. The only
non-trivial algorithm known for this case is the standard Fourier
based algorithm, described in Section~\ref{sec:prelims}. The analysis
of the running time of this algorithm reduces to the following
question:
\begin{quote}
What is the smallest $t$ such that every symmetric boolean function on $k$ variables,
which is not a constant or a parity function,
has a non-zero Fourier coefficient of order at least $1$ and at most $t$?
\end{quote}

A bound of $t_0$ implies a running time of roughly $n^{t_0}$.
A bound of $\tfrac{2k}{3}$ was provided in~\cite{MOS}.
This was improved to $\tfrac{3k}{31}$ in~\cite{LMMV}.
Here we show a bound of $O(k / \log k)$ (Theorem \ref{th:main}),
giving the first algorithm for learning symmetric $k$-juntas in time $n^{o(k)}$.

\subsection*{Techniques}
Our techniques involve a mix of number theory, combinatorics and
probability.  We start by reducing our problem to finding 0/1
solutions to a system of Diophantine equations involving binomial
coefficients, as in~\cite{LMMV}. We then take a departure
from~\cite{LMMV} by further reducing this to the problem of showing
that a certain integer-valued polynomial $P$ is constant over the set
$\{0,1,...,k\}$. We manage to prove this in two steps: First, we show
that $P$ is constant over the union of two small intervals
$\{0,...,t\}\cup\{k-t,...,k\}$. This is obtained by looking at $P$
modulo carefully chosen prime numbers. To choose these prime numbers
we use the Siegel-Walfisz theorem on the density of primes in
arithmetic progressions with modulus of moderate growth. In the second
step, we extend the constant nature of $P$ to the whole interval
$\{0,...,k\}$ by repeated applications of Lucas' Theorem. One
additional interesting aspect of our proof is the use of an
equivalence between a) the vanishing of Fourier coefficients and b)
the equality of moments of certain random variables under the uniform
measure on the hypercube and under the measure defined by the function
itself. This equivalence helps us eliminate a lot
of case analysis.

\section{Preliminaries}
\label{sec:prelims}
\subsection*{Symmetric Juntas}
Given a boolean function $f$ on $n$ variables $x_1,...,x_n$, we will
say that $x_i$ is a {\it relevant} variable for $f$ if there exist $x,y \in
\{0,1\}^n$ which differ only in the $i$-th coordinate and $f(x)\neq
f(y)$. Variables that are not relevant are called irrelevant. We will call $f$ a $k$-junta if $f$ has at most $k$ relevant
variables.

We consider the class of symmetric
juntas. A boolean function $f:\{0,1\}^k \rightarrow \{0,1\}$ on $k$
variables is a symmetric function if for any permutation $\pi\in S_k$,
$f(x_1,...,x_k) = f(\pi(x_1),...,\pi(x_k))$. Hence the value of $f$ at
$(x_1,...,x_k)$ depends only on the {\it weight} of $(x_1,...,x_k)$,
which is the number of variables that are set to $1$. A symmetric
$k$-junta is a function on $n$ variables which is symmetric on the $k$
variables it depends on.

We will describe
a symmetric boolean function on $k$ variables by a $(k+1)$-bit string
$f_0f_1...f_k$, where $f_i$ is the value of $f$ on an input of weight
$i$. The following four special symmetric functions on $k$ variables
will appear often: the two constant functions $\zero$ and $\one,$ the
parity function $\pty,$ and its complement $\ptyc$. 

\subsection*{Learning in the PAC model}
We consider the PAC learning model~\cite{Valiant}, in which we wish to learn a {\em Concept Class} $\mathcal{C} = \bigcup _{n}
\mathcal{C}_n,$ where each $\mathcal{C}_n$ is a collection of boolean
functions from $\{0,1\}^n \to \{0,1\}$. In our case, $\mathcal{C}_n$ is the class of symmetric $k$-juntas on $n$ variables. Let $\epsilon$ be an {\em accuracy
parameter} and $\delta$ a {\em confidence parameter}. A learning
algorithm $\mathcal{A}$ for $\mathcal{C}$ has access to an {\em
oracle} for $f \in \mathcal{C}_n$. A query to
the oracle outputs a labeled example $\langle x,f(x) \rangle,$
where $x$ is drawn from $\{0,1\}^n$ according to some probability
distribution $\mathcal{D}$. $\mathcal{A}$ is said to be a learning algorithm for the
class $\mathcal{C}$ under the distribution $\mathcal{D}$ if for all $f \in \mathcal{C}$, it outputs, with probability at
least $1-\delta$, a hypothesis $h$ such that $\mathrm{Pr}_{x}
[h(x)=f(x)] \geq 1-\epsilon$. We will be concerned
only with the uniform distribution and we will obtain an algorithm with accuracy parameter $\epsilon = 0$, i.e., we identify the exact function $f$.

\subsection*{Fourier Transform}
We will consider functions of the form: $f:\{0,1\}^n \rightarrow
\mathbb{R}$. 
An orthonormal basis for the functions defined on the
Boolean cube can be given by the {\it characters} of the group
$Z_2^n$. In particular, for every $S\subseteq \{1,...,n\}$, define the
following function:
$$ \chi_S(x) = (-1)^{\sum_{i\in S} x_i}.$$ 

Any real-valued function on the Boolean cube can be expressed as a
linear combination of the functions $\chi_S$.
Given $f$, we have that
$f(x) = \sum_S \hat{f}(S) \chi_S(x)$, where $\hat{f}(S)$ is the
Fourier coefficient of $f$ at $S$ and is equal to the inner product of
$f$ with $\chi_S$:
$$ \hat{f}(S) = \frac{1}{2^n}\sum_{x\in\{0,1\}^n} f(x)\chi_S(x).$$

\subsection*{Fourier-based Learning}
Let $f$ be a $k$-junta. It is known that we can exactly calculate the Fourier coefficients of
$f$ in the uniform distribution PAC model, with confidence
$1-\delta$ in time $poly(2^k,n,\log{\tfrac{1}{\delta}})$, using standard Chernoff-Hoeffding bounds
(see~\cite{LMN,MOS}). Observe further, that if $x_i$ is an irrelevant
variable for a $k$-junta $f$, then for any $S\subseteq \{x_1,...,x_n\}$ containing
$x_i$, $\hat{f}(S) = 0$. Hence if $\hat{f}(S)\neq 0$, for some $S$,
then $S$ contains only relevant variables. 

This suggests the following algorithm: Starting with $l=1$, compute
the Fourier coefficients of all subsets of $\{x_1,...,x_n\}$ of size
$l$. Collect the union of all relevant variables that correspond to
subsets with non-zero Fourier coefficients. Stop as soon as you
collect all $k$ relevant variables.

Since the function is symmetric, for any two sets $S, T$ of relevant
variables such that $|S| = |T|$, we have $\hat{f}(S) =
\hat{f}(T)$. Hence the first time that we will identify some relevant
variables in the algorithm, we will actually be able to identify all
the relevant variables. Once we find the relevant variables, finding
the truth-table of the function can be done in time $poly(2^k, \log{\tfrac{1}{\delta}})$.

The above algorithm would take time roughly $n^k$ for
$f\in\{\zero,\one,\pty,\ptyc\}$. However, these particular functions
are well known to be learnable in time $poly(n,\log{\tfrac{1}{\delta}})$. Hence the following
is true:

\begin{fact}
If every symmetric function $f\not\in\{\zero,\one,\pty,\ptyc\}$ has
a non-zero Fourier coefficient of order between 1 and $t$, then we can
learn symmetric $k$-juntas in time $n^t~poly(2^k,n,\log{\tfrac{1}{\delta}})$.
\end{fact}

\section{Main Section}

\subsection{An Equivalent Formulation}
\label{subsec:equiv}
We state an equivalent condition for the existence of a non-zero
Fourier coefficient of a boolean function $f$, as proved in~\cite{LMMV}. 
Let $f: \{0,1\}^k \to \{0,1\}$ be a
boolean function. For a vector
${\bf x}=(x_1,\ldots,x_k),$ and a set $S \subseteq [k]$, let ${\bf x}_S$ be
the projection of ${\bf x}$ on the indices of $S$. Let $\sigma \in
\{0,1\}^{|S|}.$ Define the following probabilities:
$$ p_{S,\sigma}(f):= \mathrm{Pr}\left[f({\bf x})=1~|~{\bf x}_S=\sigma  \right]$$

Unless mentioned, all
probabilities are over the uniform distribution on $\{0,1\}^k$. For $t\geq 1$, call a boolean function $f$ on $k$ variables  $t$-{\em null}, if for all sets $S \subseteq [k],$ with $ |S| = t,$ and for all $\sigma \in \{0,1\}^t,$ the probabilities
$p_{S,\sigma}(f)$ are all equal to each other.
The following lemma reveals the connection with the Fourier coefficients of $f$.
\begin{lemma}\label{lemma:nullfourier}
\cite{LMMV}~Let $f$ be a  boolean function on $k$ variables.
Then $f$ is $t$-null
for some $1\leq t \leq k,$ if and only if, for all $\emptyset \neq S
\subseteq [k]$ with cardinality at most $t$, $\hat{f}(S)=0.$
\end{lemma}
It is clear that if $s\le t$ and $f$ is $t$-null then it is also $s$-null.
\eat{
 For some $\sigma \in
\zt$, and  $T\subseteq [k], |T| = t$, let $f_{S,\sigma}$ be the function
obtained by substituting $\sigma$ for $\bf x_T$ in $f$. 
By chain rule for conditional probabilities,
$$\ptf = \frac{\mathrm{Pr}\left[f({\bf x})=1|{\bf x}_T = \sigma
\right]\cdot\mathrm{Pr}\left[{\bf x}_T=\sigma\right]}{\mathrm{Pr}\left[f({\bf x})=1
\right]} = \frac{\mathrm{Pr}\left[g({\bf x})=1\right] \cdot 2^{-t}}
{\mathrm{Pr}\left[f({\bf x})=1 \right]}$$ Hence, $\ptf = 2^{-t}$ if and
only if $\mathrm{Pr}\left[g({\bf x})=1\right]  = \mathrm{Pr}\left[f({\bf x})=1
\right]$, or equivalently, $\ff=\gf$. This is true if and only if
$\fs =0$ for all $\emptyset \neq S\subseteq T$. The lemma follows. 
}


\eat{
\noindent
\smallskip
The following  is an immediate corollary of this lemma.
\begin{corollary}\label{lemma:stnull}
Let $f$ be a  boolean function on $k$ variables. If $f$ is $t$-null for some $1\leq t \leq n$ then $f$ is $s$-null for $1\leq s \leq t.$
\end{corollary}
}
When we consider the case of symmetric functions, $p_{S,\sigma}(f)$ just depends on $t:=|S|$ and the weight $w$ of $\sigma$. 
We denote this by $p_{t,w}(f).$ It is clear that:
\begin{equation}
\label{eq:probs}
p_{t,w}(f) =
\frac{1}{2^{k-t}} \sum_{i=0}^{k} f_i{k-t\choose i - w}
\end{equation}
\eat{For a symmetric function $f$, $\mathrm{Pr}\left[f({\bf x})=1\right] =
2^{-k} \sum_{i=0}^k\bki f_i.$ Further, if ${\bf x}_s$ denotes ${\bf x}_{\{1,\ldots,s\}}$ and $\sigma$ is a string of weight $w,$ it is easy to see that $$ \mathrm{Pr}\left[f({\bf x})=1
\wedge {\bf x}_s = \sigma \right] =2^{-k} \sum_{i=0}^k{k-s \choose i-w
}f_i.$$}
where $l \choose m$ is $0$ if $m<0$ or $m>l$, and $0 \choose 0$ is $1$.
It follows that $f$ is $t$-null if for $0\leq w\leq t$, $p_{t,w}(f)$ are all equal. 
It is easy to see that the constant boolean functions $\{\zero,\one\}$ are $t$-null for all $t$ with $1\leq t\leq k$. The parity functions $\{\pty,\ptyc\}$ are also $t$-null for all $t$ satisfying $1\leq t< k$.
From Lemma~\ref{lemma:nullfourier} and Equation~\ref{eq:probs} we get: 
\begin{corollary}
\label{cor:fourier-equations}
All symmetric boolean functions $f \not\in  \{\zero,\one,\pty,\ptyc\}$ have
a non-zero Fourier coefficient of order at most $t_0$ (and at least $1$)
iff 
$\{\zero,\one,\pty,\ptyc\}$ are the only solutions to
\begin{equation}
\label{window}
\sum_{i=0}^{k-t_0} f_i{k-t_0\choose i} =
\sum_{i=1}^{k-t_0+1} f_i {k-t_0\choose i-1} =
\cdots = \sum_{i=t_0}^{k} f_i{k-t_0\choose i-t_0}\\
\end{equation}
\end{corollary}
In the next section, we show that this is true for $t_0 \le C k \bigl / \log k$ for large enough $k$. 

\vspace{.2in}
\subsection{A bound of $O(k \bigl / \log k)$.}

The following is our main theorem.
\begin{theorem}\label{th:main}
There is an absolute constant $C>0$ such that for large $k$, 
every symmetric boolean function $f$ on $k$ bits with $f \not\in  \{\zero,\one,\pty,\ptyc\}$
has a non-zero Fourier coefficient of order at most $C k \bigl/ \log k$ and at least $1$.
\end{theorem}
The rest of this section is devoted to proving Theorem~\ref{th:main}.
Suppose $f$ is a boolean function on $G=\ZZ_2^k$, such that all its Fourier
coefficients of order up to $k-N$ are $0$.
Then the values $f_j$ of $f$ satisfy \eqref{window} with $t_0 = k-N$, which, changing parameters,
can be rewritten as:
\beql{start}
\sum_j {N \choose j} f_{\nu+j} = c_N,\ \ \mbox{for all $\nu=0,\ldots,k-N$}.
\eeq 
We want to show that if $k-N \ge C k \bigl / \log k$,
for some appropriately large constant $C>0$,
then $f_j$ is either constant or alternates between $0$ and $1$.
We prove this for all $k$ sufficiently large.

Define $X_j = f_{j+1}-f_j$, for $j=0,\ldots,k-1$, and observe that the sequence $X_j$ satisfies
the homogeneous version of \eqref{start}:
\beql{start-x}
\sum_j {N \choose j} X_{\nu+j} = 0,\ \ \mbox{for all $\nu=0,\ldots,k-N-1$}.
\eeq

\noindent{\bf Remark.} In \eqref{start-x} the number $N$ can be replaced by any other integer $N_1$
in the interval $[N,k]$. This follows since all the non-constant Fourier coefficients up to order $k-N$ are
$0$.

From \eqref{start-x} the sequence $X_j$ may be defined for all $j\in\ZZ$ and $X_j \in \ZZ$
for all $j$.
From the theory of recurrence relations we know then that the sequence $X_j$ may be written
as a linear combination of the following sequences:
$$
(-1)^j, (-1)^j j, (-1)^j j^2, \ldots, (-1)^j j^{N-1}.
$$
The reason for this is that $-1$ is the only root of the characteristic polynomial
of the recurrence,
$\phi(z) = \sum_j {N \choose j} z^j = (1+z)^N$.
Therefore there is a polynomial $P(x)$, of degree at most $N-1$, such that
$$
X_j = (-1)^j P(j),\ \ \mbox{for all $j\in\ZZ$}.
$$
Clearly $P(x)$ takes integer values on integers and in particular $P(j) \in \Set{-1, 0, 1}$ for
$j=0,\ldots,k-1$.
From the well known characterization of integer-valued polynomials it follows that we may write
\beql{integer-valued}
P(x) = \sum_{j=0}^{N-1} a_j {x \choose j},\ \ \mbox{with $a_j \in \ZZ$}.
\eeq
If $p\ge N$ is a prime, and since all the factors that appear in denominators in \eqref{integer-valued}
are strictly less than $p$ (hence invertible mod $p$),
it follows that the sequence $P(j) \bmod p$, $j\in\ZZ$, may be viewed
as a polynomial with coefficients in $\ZZ_p$ and therefore is a $p$-periodic sequence mod $p$, i.e.
\beql{periodicity-modular}
P(j+p) = P(j) \bmod p,\ \ \mbox{for all $j\in\ZZ$ and $p\ge N$}.
\eeq
If, in addition, $0\le j < j+p <k$, when all $P$-values that appear in \eqref{periodicity-modular}
are in $\Set{-1,0,1}$, it follows that we have the non-modular equality
\beql{periodicity}
P(j+p) = P(j),\ \ (N \le p \le j+p <k).
\eeq
We want to show that $f\in  \{\zero,\one,\pty,\ptyc\}$. Since $X_j = f_{j+1} - f_j$ it is enough to show that either $X_j$ is identically $0$ or that
$X_j = (-1)^j$ or $X_j = (-1)^{j+1}$.
This is equivalent to showing that $P$ is a constant polynomial, constantly equal to $-1, 0$ or $1$.\\

\noindent
{\bf Notation.}\\
1. In what follows we repeatedly use the letter $C$ to denote a positive constant which depends
on no parameter (unless we say otherwise). As is customary,
this constant $C$ need not be the same in all its occurences.\\
2. We define $\epsilon$ by the relation $k-N = \epsilon k$ and assume $\epsilon \ge C \bigl / \log k$, with
$C$ a large enough positive constant.\\

We shall need various primes in intervals from now on.
The version of the prime number theorem that we will be using is
the Siegel-Walfisz theorem (see \cite[Theorem 2]{kumchev}).
Define the logarithmic integral
$$
\Li{x} = \int_2^x \frac{dt}{\log t} \sim \frac{x}{\log x},\ \ (x \to \infty).
$$
The Euler function $\varphi(q)$ denotes the number of moduli mod $q$ which are coprime to $q$.
\begin{othm} {\bf (Siegel-Walfisz)} \label{th:siegel-walfisz}
Let $\pi(x;M,a)$ be the number of primes $\le x$ which are equal to $a \bmod M$ and
assume that $(M,a)=1$.
Then if $M \le (\log x)^A$, $A$ a constant, we have
\beql{prime-asymptotics}
\pi(x;M,a) = \frac{\Li{x}}{\varphi(M)} + O(x \exp(-c\sqrt{\log x}),\ \ \mbox{(as $x \to \infty$)}.
\eeq
where $c$ depends on $A$ only (the constant in the $O(\cdot)$ term is absolute).
\end{othm}
For $\pi(x)$, the number of primes up to $x$ without any restriction,
the prime number theorem says $\pi(x) = \Li(x) + O(x \exp(-c\sqrt{\log x})$,
for some constant $c$.

These theorems guarantee that, for $x \to \infty$, the interval $[x,x+\Delta]$ has the
``expected'' number of primes whenever $\Delta \ge C x \bigl / (\log x)^A$,
whatever the constant $A$, even if we impose the condition that these primes
are equal to $a \bmod M$, as long as $M \le (\log x)^B$, for any constant $B$.

We use the above theorems along with the $p$-periodicity of $P$ to deduce that $P$ is in fact $2$-periodic on the union of $2$ small sub-intervals of $[0,k-1]$. 
\begin{lemma}\label{lm:periodic}
The polynomial $P$ satisfies the 2-periodicity condition
$$
P(j) = P(j+2),
$$
whenever $j, j+2 \in {\mathcal A} = [0,k-N] \cup [N,k-1]$.
\end{lemma}
\begin{proof}
Assume $q<r$ are two primes in $[N, N+h]$,
where $h=(k-N)/3 = \frac{\epsilon}{3} k$.
(The length of the interval $[N, N+h]$ is large enough for the prime number
theorem to guarantee the existence of many primes in it.)
From \eqref{periodicity} it follows that the finite sequences
$$
P(0),\ldots, P(k-q)\ \ \mbox{and}\ \ P(q), \ldots, P(k)
$$
are identical.
Applying \eqref{periodicity} again with $r$ we get that the finite sequences
$$
P(0),\ldots, P(k-r)\ \ \mbox{and}\ \ P(r), \ldots, P(k)
$$
are identical.
It follows that
\beql{diff-periodicity}
P(j+r-q) = P(j),\ \ \mbox{if $N+h\le j \le N+2h$ and $r>q$ primes in $[N,N+h]$}.
\eeq
We now assume that the difference $M=r-q$ is the smallest difference between two primes
in $[N,N+h]$. By the prime number theorem $M \le C \log k$.
Hence, we can apply Theorem \ref{th:siegel-walfisz}.
Since $\varphi(M) \le M \le C \log k$ in that case
Theorem \ref{th:siegel-walfisz} guarantees that the number of primes equal to $a \bmod M$
in $[N,N+h]$ is at least
$$
C\frac{h}{\log^2 k} \sim C \frac{k}{\log^3 k},
$$
whenever $(M,a)=1$. All that matters here is that this number is positive.

Let $t \in [N,N+h]$ be the smallest prime which is equal to $-1 \bmod M$.
By Theorem \ref{th:siegel-walfisz}, applied to $M$ and $-1$, its existence is guaranteed and furthermore that
$t \sim N$.
The same theorem guarantees that we can find a prime $s \in (t, N+h]$ such that
$s = 1 \bmod M$.
Then $s-t = 2 \bmod M$ or $s-t = \ell M + 2$, for some nonnegative integer $\ell$.
Therefore, for $N+h \le j \le N+2h$ we have
\begin{eqnarray*}
P(j)
	&=& P(j+s-t)\ \ \mbox{(applying \eqref{diff-periodicity} for the primes $s,t$)}\\
	&=& P(j+\ell M +2)\\
	&=& P(j+(\ell-1)M+2)\ \ \mbox{(applying \eqref{diff-periodicity} for the primes $r,q$)}\\
	&\ & \cdots\\
	&=& P(j+2).
\end{eqnarray*}
This $2$-periodicity
\beql{2-periodicity}
P(j)=P(j+2)
\eeq
is transferred to
all $j,j+2 \in {\mathcal A}$ by using \eqref{periodicity} repeatedly for appropriate primes $p$.
\end{proof}

Notice that in the sequence $X_j$, if one erases the $0$'s then one sees an alternation
of $-1$ and $1$ (this follows from the fact that $f_j \in \Set{0,1}$).
This property greatly reduces the number of allowed patterns in $X_j$ and in fact it implies that $P$ is constant in ${\mathcal A}$.
\begin{lemma}\label{lm:constant}
The polynomial $P$ is constant in ${\mathcal A}$ (defined in Lemma \ref{lm:periodic}).
\end{lemma}
\begin{proof}
From Lemma \ref{lm:periodic} the values of $P$ in $[N,k-1]$ must be a $2$-periodic sequence.
The only essentially different non-constant $2$-periodic patterns for the values of $P$ in $[N,k-1]$ are
$010101\ldots$ and $(-1)1(-1)1\ldots$ and they both violate the property that $X_j=(-1)^j P(j)$ must satisfy,
namely that if one erases the $0$'s then one must see an alternation of $1$ and $-1$.
Therefore $P$ is constant in each of the two intervals of ${\mathcal A}$.
From the $p$-periodicity \eqref{periodicity} it follows
that the constant is the same in both intervals.
\end{proof}

We now extend the set on which $P$ is constant to a superset of ${\mathcal A}$ that contains a small interval around $k/2$. 
We will make use of the following theorem which follows from Lucas' Theorem \cite[Ch.\ 3]{cameron}.
\begin{theorem}\label{th:lucas}
If $r$ is a prime which does not divide $n$ then ${m r \choose n} = 0 \bmod r$.
Also, if $0\le m < r$ then ${mr \choose lr} = {m \choose l} \bmod r$.
\end{theorem}

\begin{lemma}\label{lm:middle}
Let $a = (1/2-\epsilon/2)k$ and $b = (1/2+\epsilon/2)k$.
Then $P(l)=P(0)$ for $a\le l \le b$.
\end{lemma}
\begin{proof}
We shall apply Theorem \ref{th:lucas} with $m=2$ and with a prime $r$ such that
$2r-N$ takes the minimal possible nonnegative value.
It follows from the prime number theorem
that $2r-N=o(\epsilon k)$.
And it follows from the remark after \eqref{start-x} that
$$
\sum_j (-1)^j {2r \choose j} P(j+\nu) = 0,\ \ \ (\nu \in \ZZ).
$$
Taking residues mod $r$ and using Theorem \ref{th:lucas} for $m=2$ we obtain
$$
P(\nu) - 2 P(\nu+r) + P(\nu+2r) = 0 \bmod r,\ \ \ (\nu \in \ZZ).
$$
By our particular choice of $r$ we have $P(\nu) = P(\nu+2r) = P(0)$ whenever
$\nu \in [0,k-N-o(\epsilon k)]$.
It follows that $P(\nu+r) = P(0)$.
Applying this for all $\nu\in[0,k-N-o(\epsilon k)]$ we get $P(l) = P(0)$
for all $l$ in the interval $(a+o(\epsilon k), b-o(\epsilon k))$.
To get rid of the $o(\epsilon k)$ terms in the interval above, just choose
a slightly larger $r$ and apply again for all $\nu \in [0,k-N-o(\epsilon k)]$.
\end{proof}

So far we have proved $P(l)=P(0)$ on the set $$ {\mathcal A}_2 =
[0,k-N] \cup [a,b] \cup [N,k-1], $$ which consists of three equispaced
intervals of roughly equal size $\epsilon k$. We consider $2$ cases for $P$. The first is when $P$ is $0$ on ${\mathcal A}_2$ and the second is when $P$ is $1$ or $-1$.

In the case that $P$ is $0$ on ${\mathcal A}_2$, we shall need the following theorem, which already gives a lot of significant information
about the function $f$.
It should be thought of as analogous to the fact that the moments of a (vector) random variable
can be read off the Fourier Transform of its distribution (the {\em characteristic function})
by looking at derivatives at $0$.
\begin{theorem}\label{th:moments}
Suppose $f: G=\ZZ_2^k = \Set{0,1}^k \to \RR$ is nonnegative (and not
identically $0$) and has all its
Fourier coefficients of order at most $r$ (and at least 1) equal to $0$.
Let $\mu$ denote the uniform probability measure on the cube $G$ and
$\nu$ denote the probability measure on $G$ defined by
$$
\nu(A) = \sum_{x \in A} f(x) \Bigl / \sum_{x \in G} f(x),\ \ \mbox{($A \subseteq G$)}.
$$
Let also $X_1,\ldots,X_k$ denote the coordinate functions on $G$, which we
view as random variables.
Then for all $i_1<i_2<\cdots<i_s$, $0\le s \le r$, we have
$$
\Mean_\nu(X_{i_1}\cdots X_{i_s}) = \Mean_\mu(X_{i_1}\cdots X_{i_s}).
$$
\end{theorem}
\begin{proof}
Let $F = \sum_{x\in G}f(x)$.
We assume for simplicity that $i_1=1,\ldots,i_s=s$.
Then, writing $x=(x_1,x_2,\ldots,x_k)$ and $[s]=\Set{1,\ldots,s}$, we have
\begin{eqnarray*}
\Mean_\nu(X_1\cdots X_s)
	&=& \frac{1}{F} \sum_{x\in G} f(x) x_1\cdots x_s \\
	&=& \frac{1}{F} \sum_{x\in G} f(x) \frac{1+(-1)^{x_1+1}}{2}\cdots\frac{1+(-1)^{x_s+1}}{2}\\
	&=& \frac{1}{2^s F} \sum_{x\in G} f(x) \sum_{S \subseteq [s]}(-1)^{\Abs{S}+\sum_{i\in S}x_i}\\
	&=& \frac{\Abs{G}}{2^s F} \sum_{S\subseteq [s]} (-1)^{\Abs{S}} \frac{1}{\Abs{G}}
		\sum_{x\in G} f(x)(-1)^{\sum_{i\in S}x_i}\\
	&=& \frac{\Abs{G}}{2^s F} \sum_{S\subseteq [s]} (-1)^{\Abs{S}} \ft{f}(S)\\
	&=& \frac{\Abs{G}}{2^s F} \ft{f}(0)\ \ \mbox{(by the vanishing of $\ft{f}$)}\\
	&=& 2^{-s}\\
	&=& \Mean_\mu(X_1\cdots X_s)
\end{eqnarray*}
\end{proof}

\noindent
{\bf Remarks.}\\
1. For functions $f:\{0,1\}^k\to\{0,1\}$,
the above theorem follows directly from the definition of $t$-nullity in Section~\ref{subsec:equiv}.
However, as we shall see in the proof of Lemma~\ref{lm:nonzero} we need to apply this theorem for functions whose range is not $\{0,1\}$.\\
2. If the nonnegative function $f$ is symmetric then the
identity of moments up to order $r$ with those of the uniform distribution ($r$-wise independence)
and the vanishing of the non-constant Fourier coefficiens of weight up to $r$ are equivalent.
This can be proved by induction on $r$. We do not use this here.

\begin{corollary}\label{cor:moments}
Under the assumptions and definitions of Theorem \ref{th:moments}
the random variable $S=X_1+\cdots+X_k$ has the same power
moments under the probability measures $\mu$ and $\nu$,
up to order $r$.
\end{corollary}
\begin{proof}
The power $S^s$, $s\le r$,
can be written as a sum of terms of the type $X_{i_1}\cdots X_{i_t}$, for $t\le s$.
One uses the fact that $X_j^2 = X_j$.
\end{proof}
\begin{lemma}\label{lm:nonzero}
If $P$ is $0$ on ${\mathcal A}_2$, then $f\in\{\zero,\one\}$.
\end{lemma}
\begin{proof}
Suppose the polynomial $P$ is constantly equal to $0$ on the set ${\mathcal A}_2$
and that $f\not\in\{\zero,\one\}$. 
The sequence $f_j$ is constant in each of the three intervals of
${\mathcal A}_2$.
By possibly considering $1-f$ (whose Fourier coefficients vanish exactly where those of $f$ do), we
may assume that $f_j=0$ on the middle interval $(a,b)$.
Define the nonnegative function $g:G \to \RR$ by
$$
g(x_1,\ldots,x_k) = f(x_1,\ldots,x_k) + f(1-x_1, \ldots, 1-x_k),
$$
and observe that the Fourier coefficients of $g$ of weight at most $k-N$ vanish.
Let $\tau$ be the distribution of the random variable $S=X_1+\cdots+X_k$ under
the measure induced by $g$ on $G$ (each vertex $x\in G$ has probability proportional to $g(x)$). Note that this is a well defined probability distribution since we assumed that $f$ and $1-f$ are not the $\zero$ function.
Clearly $\tau$ is symmetric about $k/2$ and has no mass in $(a,b)$,
since both $f(x_1,\ldots,x_k)$ and $f(1-x_1,\ldots,1-x_k)$ are $0$ when $x_1+\cdots+x_k \in (a,b)$.
The $s$-th moment with respect to the measure $\tau$
of the variable $S$ in Corollary \ref{cor:moments} is the expression
$$
M(\tau, s) = \frac{1}{F} \sum_j g_j {k \choose j} j^s,
$$
where again $F = \sum_j g_j {k \choose j}$.
By Corollary \ref{cor:moments} this must equal the $s$-th moment with respect to
the binomial measure $\mu$, which is the quantity
$$
M(\mu, s) = 2^{-k} \sum_j {k \choose j} j^s.
$$
But the variance of $S$ under $\mu$ is
$$
M(\mu, 2) - M(\mu, 1)^2 = k,
$$
since under $\mu$ the random variables $X_1,\ldots,X_k$ are independent,
while the variance of $S$ under $\tau$ is
$$
M(\tau, 2)-M(\tau, 1)^2 \ge C \epsilon^2 k^2,
$$
as half the mass of $\tau$
sits to the left of $\frac{1-\epsilon}{2}k$ and half to
the right of $\frac{1+\epsilon}{2}k$.
These orders of magnitude are different whenever $\epsilon \ge C \bigl / \sqrt k$, which
is true in our case as $\epsilon \ge C \bigl / \log k$.
This contradiction proves that $P$ cannot equal $0$ on ${\mathcal A}_2$.
\end{proof}

\noindent
{\bf Extending ${\mathcal A}_2$ to $[0,k-1]$.}

The rest of the proof goes as follows.
By Lemma \ref{lm:nonzero}, we may assume that
$P(l)=1$ or $-1$ for $l \in {\mathcal A}_2$. Without loss of generality, assume $P$ is $1$ on ${\mathcal A}_2$.
We apply Theorem \ref{th:lucas}
for $m=4,8,16,\ldots$ successively and each time we choose a prime $r$ such that
$mr-N$ is minimized.
Theorem \ref{th:lucas} gives for all $\nu \in \ZZ$
\beql{long-relation}
P(\nu)-mP(\nu+r)+{m \choose 2}P(\nu+2r) - \cdots +P(\nu+mr) = 0 \bmod r.
\eeq
When $\nu \in [0,k-N]$
the numbers $\nu+lr$ for even $l$ in \eqref{long-relation} are in the set ${\mathcal A}_{m/2}$ and therefore
the corresponding $P$ values are all $1$, by induction on $m$.
In order to deduce that \eqref{long-relation} holds as an identity of integers (not residue classes)
it is enough to guarantee that the sum of the absolute values of all terms is less than $r$.
This amounts to the inequality $2^m < r$. Given that $mr \sim k$ this is true if we can guarantee that
\beql{m-bound}
m \le c_1 \log k,
\eeq
for some small enough constant $c_1$.
Therefore, as long as $m$ satisfies the bound \eqref{m-bound},
we have that, for $\nu \in [0,k-N]$,
\beql{long-relation-non-modular}
P(\nu)-mP(\nu+r)+{m \choose 2}P(\nu+2r) - \cdots +P(\nu+mr) = 0.
\eeq
Since the total weights of the positive and
negative terms in \eqref{long-relation-non-modular} are the same,
it follows that the $P(\nu+lr)$ terms corresponding to odd $l$ are also
$1$.

Each time we perform this operation we deduce that $P$ is $1$ on a collection
of intervals ${\mathcal A}_m$ which consists of ${\mathcal A}_{m/2}$ and one
interval of length $\epsilon k$ in the middle of the gap between
any two succesive intervals of ${\mathcal A}_{m/2}$.
So ${\mathcal A}_m$ has $m+1$ disjoint equispaced intervals of length $\epsilon k$.
We apply this operation until we have $\epsilon m \sim 1$, which implies
that we have covered the whole interval $[0,k-1]$ with our set ${\mathcal A}_m$.
We need to make sure that \eqref{m-bound} still holds then.
Since $\epsilon m \sim 1$ this is achieved by setting
$\epsilon = C \bigl / \log k$, for a large enough constant $C$.
At the end of this process, there could still be some very small possibly uncovered intervals of size $o(\epsilon k)$. However since we have already shown that $P(l) = 1$ on a set of $k-o(\epsilon k)$ entries, we can use the fact that $P$ has degree at most $N-1$ to obtain that $P(l) = 1$ on the whole interval $[0,k-1]$. 

This concludes the proof of the Theorem~\ref{th:main}, which implies:
\begin{corollary}
The class of symmetric $k$-juntas can be learned exactly under the
uniform distribution with confidence $1-\delta$ in time
$n^{O(k/\log{k})}\cdot \mathrm{poly} (2^k,n,\log (1/\delta))$.
\end{corollary}

\section{Discussion}
The main open question is to obtain tight upper and lower bounds on the running time of the Fourier-based algorithm for symmetric juntas. It may even be that for large $k$, every symmetric function has a non-zero Fourier coefficient of constant order.

It should also be noted that in the case of balanced symmetric functions, i.e., symmetric functions with $Pr[f(x) = 1] = 1/2$, a bound of $O(k^{0.548})$ follows from~\cite{GR} (see~\cite{MOS}). Hence to improve our result, one may focus on finding new techniques for unbalanced functions.
\bibliographystyle{abbrv}
\bibliography{boolean}
\end{document}